\documentclass[12pt,twoside,reqno]{amsart}
\usepackage{pifont}
\usepackage{amsmath}
\usepackage{amsthm}
\usepackage{bm}
\usepackage{amsfonts}
\usepackage{amssymb}
\usepackage{enumitem}
\usepackage{caption}   
\usepackage{booktabs}  
\usepackage{latexsym}
\usepackage{amstext}
\usepackage{array}
\usepackage{graphicx}
\usepackage{tikz}
\usepackage{mathrsfs}
\usepackage{mathtools}

\date{}
\allowdisplaybreaks[4] \footskip=15pt
\renewcommand{\uppercasenonmath}[1]{}

\newtheorem{thm}[subsection]{Theorem\;}
\newtheorem{cor}[subsection]{Corollary\;}
\newtheorem{Def}[subsection]{Definition\;}
\newtheorem{lem}[subsection]{Lemma\;}
\newtheorem{remark}{Remark\;}

\newtheorem{prop}[subsection]{Proposition\;}
\newtheorem{exm}[subsection]{Example\;}
  
\newcommand{\bthm}{\begin{thm} }
	\newcommand{\ethm}{\end{thm} }
\newcommand{\bpro}{\begin{prop}}
	\newcommand{\epro}{\end{prop}}
\newcommand{\modclass}[1]{\ (\mathrm{mod}\ #1)}
\newcommand{\bdf}{\begin{Def}}
	\newcommand{\edf}{\end{Def}}
\newcommand{\bexm}{\begin{exm}}
	\newcommand{\eexm}{\end{exm}}

\newcommand{\blem}{\begin{lem}}
	\newcommand{\elem}{\end{lem}}
\newcommand{\bpf}{\begin{proof}}
	\newcommand{\epf}{\end{proof}}
\newcommand{\bcor}{\begin{cor}}
	\newcommand{\ecor}{\end{cor}}
\newcommand{\ba}{\begin{array}}
	\newcommand{\ea}{\end{array}}
\newcommand{\bea}{\begin{eqnarray}}
	\newcommand{\eea}{\end{eqnarray}}

\newcommand{\brem}{\begin{remark}}
	\newcommand{\erem}{\end{remark}}

\graphicspath{dir-list}
\begin{document}
\begin{center}
{\large \bf On integers of the form \(p+F_{2^k}+F_q\)}
\vskip 0.8cm
{\small Yang Gao}\\
{\small School of Mathematics and Statistics, Huaiyin Normal University}\\
{\small Huaian, Jiangsu 223300, P.R. China}\\
{\small Email: \texttt{dg21210003@smail.nju.edu.cn}}
\end{center}
{\bf Abstract }
In 1934, Romanoff proved that the set of positive integers representable 
as the sum of a prime and a power of two has positive lower density.
Erd\H{o}s later constructed an infinite arithmetic progression of odd integers
none of which admits such a representation.
Let \(F_n\) be the Fibonacci sequence. In this paper, we prove that the set of integers of the form \(p+F_{2^k}+F_q\), where \(p,q\)
are primes and \(k\ge0\), has positive lower asymptotic density.
The same holds for the set of integers not of this form.

{\small 
{\bf Mathematics Subject Classification 2020.} 11P32, 11A41, 11B13.}

{\bf Keywords.} Romanoff type problems; primes; asymptotic density; application of the sieve methods.
\section{Introduction and main results}\label{02}
Let $\mathcal{P}$ be the set of all primes and $\mathbb{N}_0=\{0,1,2,\dots\}$.
In 1849, de Polignac~\cite{dePolignac1849} conjectured that every odd integer
greater than $3$ can be expressed as the sum of a prime and a power of two,
but soon discovered a counterexample.  Despite this, Romanoff~\cite{Romanoff1934}
proved in 1934 that the set of odd integers representable in this way has
positive lower asymptotic density.  To study the complementary set,
Erd\H{o}s~\cite{Erdos1950} introduced the notion of a covering system of
congruences and used it to construct an infinite arithmetic progression of
odd integers that cannot be written as $p+2^n$, thereby showing that the
complement also has positive lower density.  These classical results gave
rise to a broad theory of Romanoff type problems; for quantitative refinements
of Romanoff's constant, local density in arithmetic progressions, and sums
with linear recurrences, see
\cite{ChenSun2004,Dubickas2013,ElsholtzSchlage2018,Sun2010,SunFang2008}.

A fruitful direction replaces the powers of two by linear recurrences.
For Fibonacci numbers $F_n$ ($F_0=0$, $F_1=F_2=1$, $F_{n+2}=F_{n+1}+F_n$),
Lee~\cite{Lee2010} first established that $p+F_k$ has positive lower density.
This was later improved by Liu and Xue~\cite{LiuXue2021} and by Wang and
Chen~\cite{WangChen2023}, the latter also proving uniqueness of representations
for a positive proportion of integers.  
 Wang~\cite{Wang2025} extended the
analogue to Lucas numbers, showing that the set $p+L_n$ also has positive
lower density. 
In the opposite direction, \v{S}iurys~\cite{Siurys2016} showed that there exists
an infinite arithmetic progression of integers not of the form \(p^\alpha\pm F_n\),
where \(p\) is prime and \(\alpha\ge 1\).  More recently, Wang~\cite{Wang2025b} used covering systems to construct an
infinite arithmetic progression consisting of integers that can be expressed
neither as \(p+F_m\) nor as \(q+L_n\).
For sums involving two or more sparse terms, early results include
Crocker~\cite{Crocker1971}, who proved that infinitely many odd integers
cannot be expressed as \(p+2^a+2^b\), Sun and Li~\cite{SunLi2001}, who
studied integers not of the form \(c(2^a+2^b)+p^\alpha\), and
Yuan~\cite{Yuan2004}, who treated similar problems.
Subsequently, Pan~\cite{Pan2011} proved that the integers not of the form
\(p+2^a+2^b\) have lower density at least \(x^{1-\varepsilon}\).
Elsholtz, Luca and Planitzer~\cite{ElsholtzLucaPlanitzer2018} then established
positive lower density for the mixed forms \(p+2^{2^k}+m!\).  More recently, Ding~\cite{Ding2022} obtained a positive
density result for \(p+2^{a^2}+2^{b^2}\), and Chen and
Xu~\cite{ChenXu2025a} sharpened this by showing that a positive proportion
of integers admit a unique representation in that form.
For sums with several sparse terms, Chen and Wang~\cite{ChenWang2024}
studied the sum of a prime and a term from certain exponential sequences.
Chen and Xu~\cite{ChenXu2024} proved that if \(r_1,\dots,r_t\) are positive
integers satisfying \(r_1^{-1}+\cdots+r_t^{-1}\ge 1\), then both the set of
integers of the form \(p+2^{k_1^{r_1}}+\cdots+2^{k_t^{r_t}}\) and the set of
integers of the form \(p+F_{k_1^{r_1}}+\cdots+F_{k_t^{r_t}}\) have positive
lower asymptotic density, where \(p\) is a prime and \(k_1,\dots,k_t\) are
nonnegative integers.
Recently, Xu and Chen~\cite{XuChen2025b} studied $n=p+F_{k_1^2}+F_{k_2^2}$
with $\{F_n\}$ the Fibonacci sequence with $F_0=0$, $F_1=1$, and proved that
both the set of integers with exactly one representation and the set with at
least two representations have positive density, while the non‑representable
set contains an infinite arithmetic progression. For further developments on Romanoff type problems in other directions,
including polynomial analogues, we refer the reader to
\cite{DaiPan2026,LiuLu2011,Sun2010}.
In this paper we consider integers of the form
\[
n = p + F_{2^k} + F_q,
\]
where $p,q\in\mathcal{P}$ and $k\in\mathbb{N}_0$.
\begin{thm}
\label{thm:neg}
There exists an infinite arithmetic progression of positive integers, none of
which can be written as \(p+F_{2^k}+F_q\), where \(p,q\) are primes and
\(k\ge 0\).  Consequently, the set of integers not representable in this form
has positive lower asymptotic density.
\end{thm}

\begin{thm}
\label{thm:pos}
The set of positive integers representable as \(p+F_{2^k}+F_q\), where
\(p,q\) are primes and \(k\ge 0\), has positive lower asymptotic density.
\end{thm}
\begin{remark}
 The analogous statements for Lucas numbers also hold, since the proof is entirely similar and relies on the same elementary properties of the sequences modulo the relevant primes.
\end{remark}
This paper is organized as follows: In Section~2, we collect some
preliminary lemmas on Fibonacci numbers and sieve estimates.
In Section~3, we prove Theorem~\ref{thm:neg} by constructing an infinite
arithmetic progression of integers not of the form \(p+F_{2^k}+F_q\).
In Section~4, we give the proof of Theorem~\ref{thm:pos}, which shows that
the representable set has positive lower density.

\section{Preliminaries}
\subsection*{2.1 The fundamental properties of the Fibonacci sequence}%
\par
The inequalities in the following lemma are elementary; their proofs are omitted.
For the congruence properties of the Fibonacci sequence modulo primes,
we refer to Somer~\cite{Somer1990}.
 
\begin{lem}\label{lem:fib}
Let $\gamma = \frac{1+\sqrt{5}}{2}$ and let $F_n$ denote the $n$‑th Fibonacci number \rm{(}$F_0=0,F_1=F_2=1$, $F_{n+2}=F_{n+1}+F_n$\rm{)}. 
Then the following statements hold.
\begin{enumerate}[label=\normalfont(\arabic*)]   
    \item $\displaystyle \left|F_n - \frac{\gamma^n}{\sqrt{5}}\right| < \frac{1}{2}$  and \(\prod\limits_{k=1}^n F_k \leq \gamma^{\frac{n(n+1)}{2}} \quad \;\text{for all}\; n \in \mathbb{N}^*\)
    \item Define $L(x) = \dfrac{\log\bigl(\sqrt{5}\,(x+\tfrac{1}{2})\bigr)}{\log\gamma},x>1$. If $F_n \leqslant x $ then $ n \leqslant L(x)$.
    \item Define $W(x) = \dfrac{\log\bigl(\sqrt{5}\,(x-\tfrac{1}{2})\bigr)}{\log\gamma},x>1$.   
           If positive integer $n \leqslant W(x)$ then $ F_n \leqslant x $.
    \item $2\log x < L(x) < 4\log x \quad (x \ge 3);\;2\log x < W(x) < 3\log x \quad (x \ge 5)
.$
    \item Every positive integer $n$ can be uniquely written as
\[
n = F_{i_1} + F_{i_2} + \cdots + F_{i_r},
\]
where $i_r \geqslant 2$ and $2 \leqslant i_j - i_{j+1}$ for all $j$; this representation is called the \emph{Zeckendorf representation} of $n$. 
    \item For any positive integer \(d\), the sequence \(F_n \bmod d\) is purely periodic: there exists a smallest positive integer \(k(d)\) (the Pisano period) such that
\[
F_{n+k(d)} \equiv F_n \pmod{d} \quad \text{for all } n \geq 0.
\] Furthermore, every prime $p$ divides the Fibonacci number $F_{k(p)}$.
    \item If $p$ is prime then in one complete period of the Fibonacci numbers taken modulo $p$, every residue occurs no more than four times.
   \item For integers \(m>n\ge 0\) such that \(m-n\) is even, the following identity holds:
\[
F_m - F_n = F_{m-1} + F_{m-3} + \cdots + F_{n+1}.
\]
\end{enumerate}
\end{lem}
We now list three tables concerning the distribution and periodicity of the
Fibonacci sequence modulo certain primes.

\begin{center}
\captionof{table}{Pisano periods $k(p)$ for the primes $p \le 50$}
\small
\renewcommand{\arraystretch}{1.3}
\begin{tabular*}{\textwidth}{@{\extracolsep{\fill}}c*{15}{c}}
\hline
$p$      & $2$ & $3$ & $5$ & $7$ & $11$ & $13$ & $17$ & $19$ & $23$ 
         & $29$ & $31$ & $37$ & $41$ & $43$ & $47$ \\
\hline
$k(p)$ & $3$ & $8$ & $20$ & $16$ & $10$ & $28$ & $36$ & $18$ & $48$ 
         & $14$ & $30$ & $76$ & $40$ & $88$ & $32$ \\
\hline
\end{tabular*}
\end{center}

\begin{center}
\captionof{table}{Fibonacci numbers $F_n \bmod 7$ for $n=0$ to $16$}
\small
\renewcommand{\arraystretch}{1.3}
\begin{tabular*}{\textwidth}{@{\extracolsep{\fill}}c*{17}{c}}
\hline
$n$           & $0$ & $1$ & $2$ & $3$ & $4$ & $5$ & $6$ & $7$ & $8$ 
              & $9$ & $10$ & $11$ & $12$ & $13$ & $14$ & $15$ & $16$ \\
\hline
$F_n \bmod 7$ & $0$ & $1$ & $1$ & $2$ & $3$ & $5$ & $1$ & $6$ & $0$ 
              & $6$ & $6$ & $5$ & $4$ & $2$ & $6$ & $1$ & $0$ \\
\hline
\end{tabular*}
\end{center}

\begin{center}
\captionof{table}{Fibonacci numbers $F_n \bmod 19$ for $n=0$ to $18$}
\small
\renewcommand{\arraystretch}{1.3}
\begin{tabular*}{\textwidth}{@{\extracolsep{\fill}}c*{19}{c}}
\hline
$n$            & $0$ & $1$ & $2$ & $3$ & $4$ & $5$ & $6$ & $7$ & $8$ 
               & $9$ & $10$ & $11$ & $12$ & $13$ & $14$ & $15$ & $16$ 
               & $17$ & $18$ \\
\hline
$F_n \bmod 19$ & $0$ & $1$ & $1$ & $2$ & $3$ & $5$ & $8$ & $13$ & $2$ 
               & $15$ & $17$ & $13$ & $11$ & $5$ & $16$ & $2$ & $18$ 
               & $1$ & $0$ \\
\hline
\end{tabular*}
\end{center}
\begin{lem}[{\cite[Theorem 328]{HardyWright1979}}]\label{lem:Hardy}
There exists an absolute constant \(c>0\) such that
\[
\varphi(x) \ge c\,\frac{x}{\log\log x}
\qquad\text{for all } x\ge 3,
\]
where \(\varphi\) is Euler's totient function.
\end{lem}

\begin{lem}[{\cite[(6.6)]{ChenXu2024}}]\label{lem:chen}
Let $k(d)$ be the minimal positive period of the Fibonacci sequence modulo $d$, and set
$v(d) = \max\{\,k(p) : p \mid d,\ p \text{ prime}\,\}$.
Then for any $\alpha > 0$, the series
\[
\sum_{d=1}^{\infty} \frac{\mu^2(d)}{d \, v(d)^\alpha}
\]
converges.
\end{lem}
\bpf To keep this paper self-contained, and since the lemma is short and its proof
can be given in a few lines, we include a short version of the proof
originally due to Chen~\cite{ChenXu2024}.  It suffices to prove
\[
\sum_{n=1}^{\infty}\frac{1}{n^\alpha}\sum_{v(d)=n}\frac{\mu^{2}(d)}{d}<\infty.
\]

For large \(N\), partial summation gives
\[
\sum_{n=1}^{N}\frac{1}{n^\alpha}\sum_{v(d)=n}\frac{\mu^{2}(d)}{d}
\le \sum_{n=1}^{N}\frac{\alpha}{n^{\alpha+1}}\sum_{v(d)\le n}\frac{\mu^{2}(d)}{d}
   +\frac{1}{(N+1)^\alpha}\sum_{v(d)\le N}\frac{\mu^{2}(d)}{d},
\tag{1}
\]
where we used \(\frac{1}{n^\alpha}-\frac{1}{(n+1)^\alpha}\le\frac{\alpha}{n^{\alpha+1}}\).

If \(v(d)\le n\), then for every prime \(p\mid d\) we have \(k(p)\le n\).
By Lemma~\ref{lem:fib}(6), \(p\mid F_{k(p)}\), so \(p\mid F_{1}\cdots F_{n}\).
 Therefore
\[
\sum_{v(d)\le n}\frac{\mu^{2}(d)}{d}
\le \prod_{p\mid F_{1}\cdots F_{n}}\Bigl(1+\frac{1}{p}\Bigr)
\le \prod_{p\mid F_{1}\cdots F_{n}}\Bigl(1-\frac{1}{p}\Bigr)^{-1}
   =\frac{F_{1}\cdots F_{n}}{\varphi(F_{1}\cdots F_{n})}.
\tag{2}
\]

Applying Lemma~\ref{lem:Hardy} to \(m=F_{1}\cdots F_{n}\) yields
the right-hand side of (2) is
\(\ll\log\log(F_{1}\cdots F_{n})\).  By Lemma~\ref{lem:fib}(1), whence
\(\log\log(F_{1}\cdots F_{n})\ll\log n\).  Substituting this into (2) gives
\[
\sum_{v(d)\le n}\frac{\mu^{2}(d)}{d}\ll\log n.
\tag{3}
\]

Inserting (3) into (1), we obtain
\[
\sum_{n=1}^{N}\frac{1}{n^\alpha}\sum_{v(d)=n}\frac{\mu^{2}(d)}{d}
\ll \sum_{n=1}^{N}\frac{\log n}{n^{\alpha+1}}+\frac{\log N}{(N+1)^\alpha}.
\]
Thus the series converges.\qed
\epf

\subsection*{2.2 Sieve estimates}
\begin{lem}[Bombieri–Davenport]\label{lem:BD}
Let $h>0$ be an even integer and let $x$ be sufficiently large. 
Let $\pi_{h}(x)$ be the number of primes $p$ such that $p\le x$ and $p+h$ is prime.
Then
\[
\pi_{h}(x)\ll\,\frac{x}{\log^2 x}\prod_{p\mid h}\Bigl(1+\frac{1}{p}\Bigr).
\]
\end{lem}
This estimate can be found in Nathanson~\cite[Theorem~7.3]{Nathanson1996}.

\begin{lem}[Brun--Titchmarsh inequality]\label{lem:BT}
Let $k$ be a positive integer and let $l$ be an integer coprime to $k$. 
If $x > k$, then the number of primes $p \le x$ with $p \equiv l \pmod k$ satisfies
\[
\pi(x; k, l) \le \frac{2x}{\varphi(k) \log(x/k)},
\]
where $\varphi$ is Euler's totient function.
\end{lem}

 \vskip 2mm
\section{Proof of Theorem \ref{thm:neg}}
\blem\label{lem:arith}
There exists an infinite arithmetic progression \(\mathcal{A}\) of positive integers such that no integer \(n\in\mathcal{A}\) can be expressed as
\[
n=p+F_{2^k}+F_q,
\]
where \(p,q\) are primes and \(k\) is a non-negative integer.
\elem

\begin{proof}
Let $
a \pmod c = \{a + ck : k \in \mathbb{Z}\}.
$
\[
\begin{aligned}
\text{Let}\;\mathcal{A}=
&0\modclass{2}\; \cap \;1\modclass{3} \;\cap \;0\modclass{5} \;\cap\; 2\modclass{7} \\
&\;\cap 3\;\modclass{11} \;\cap\; 12\modclass{19} \;\cap\; 0\modclass{23} \;\cap\; 1\modclass{47}.
\end{aligned}
\]
By the Chinese remainder theorem, \(\mathcal{A}\) is a single residue class modulo
\[
M=2\cdot 3\cdot 5\cdot 7\cdot 11\cdot 19\cdot 23\cdot 47.
\]

Suppose that
\(
n=p+F_{2^k}+F_q
\)
with $n\in \mathcal{A},$ primes \(p,q\) and integer \(k\ge 0\).
Since \(3\nmid 2^k\) and the Fibonacci sequence modulo \(2\) has period \(3\), it follows that \(F_{2^k}\) is odd for every \(k\ge 0\).
As \(n\in 0\modclass{2} \) , the parity of \(p+F_{2^k}+F_q\) forces either
\[
p\ \text{is odd and } q=3,
\qquad\text{or}\qquad
p=2\ \text{and } q\ne 3.
\]

First suppose that \(p\ge 3\) and \(q=3\). Then \(n=p+F_{2^k}+F_3=p+F_{2^k}+2\).

If \(k=\!0\) or \(1\), then \(F_{2^k}=1\), so \(n=p+3\). From \(n\equiv 3\pmod{11}\), we obtain \(p\equiv 0\pmod{11}\), hence \(p=11\). Thus \(n=14\), contradicting \(n\equiv 2\pmod 7\).

If \(k=2\), then \(F_4=3\), so \(n=p+5\). Since \(n\equiv 0\pmod 5\), we get \(p\equiv 0\pmod 5\), whence \(p=5\) and \(n=10\). But \(10\not\equiv 2\pmod 7\), a contradiction.

If \(k=3\), then \(F_4=3\), so \(n=5+F_q\). From \(n\equiv 12\pmod{19}\), we obtain \(F_q\equiv 7\pmod{19}\), but Table 3 shows that \(7\) is not a Fibonacci residue modulo \(19\), contradiction.

If \(k=\!4\), then \(F_8=21\), hence \(n=23+F_q\). From \(n\equiv 0\pmod{23}\), we get \(F_q\equiv 0\pmod{23}\). Recalling from Table 1 that the Fibonacci sequence modulo \(23\) has period \(48\), and by an elementary computation we know that for \(0\le m\le 47\), \(23\mid F_m\) if and only if \(m=0,24\), we deduce \(24\mid q\), impossible for a prime.

If \(k\ge 5\), then \(32\mid 2^k\) and \(16\mid 2^k\). From Table 1 we recall that the Fibonacci sequence has period \(16\) modulo \(7\) and period \(32\) modulo \(47\). Consequently,
\[
F_{2^k}\equiv 0\pmod 7,\qquad F_{2^k}\equiv 0\pmod {47}.
\]
Since \(n\equiv 2\pmod 7\), it follows that
\[
p+2\equiv 2\pmod 7,
\]
so \(p\equiv 0\pmod 7\), whence \(p=7\). Thus \(n=9+F_{2^k}\). Reducing this congruence modulo \(47\) and using \(F_{2^k}\equiv 0\pmod {47}\), we obtain \(n\equiv 9\pmod {47}\), which contradicts the required congruence \(n\equiv 1\pmod {47}\). Therefore the first case is impossible.

Now we only need to turn to the case where \(n=p+F_{2^k}+F_q\) with \(p=2\) and \(q\neq 3\).

If \(q=2\), then \(n=3+F_{2^k}\).  If \(k=0\) or \(1\), then \(F_{2^k}=1\) and \(n=4\), which contradicts \(n\equiv 0\pmod 5\). If \(k\ge 2\), then \(4\mid 2^k\). Recalling from Table 1 that the Fibonacci sequence modulo \(3\) has period \(8\), and that \(3\mid F_m\) for \(0\le m\le 7\) if and only if \(m=0,4\), we obtain \(F_{2^k}\equiv 0\pmod 3\) for all \(k\ge 2\). Thus \(n\equiv 0\pmod 3\), contradicting \(n\equiv 1\pmod 3\). Hence \(q\ne 2\), so \(q\ge 5\).

Now we only need to turn to the case where \(n=p+F_{2^k}+F_q\) with \(p=2\) and \(q\geq 5\).

If \(k=\!0\) or \(1\), then \(F_{2^k}=1\), so \(n=3+F_q\). From \(n\equiv 12\pmod{19}\), we get \(F_q\equiv 9\pmod{19}\). However, from Table 3 we know that \(9\) is not a Fibonacci residue modulo \(19\), a contradiction.

If \(k=2\), then \(F_4=3\), so \(n=5+F_q\). From \(n\equiv 12\pmod{19}\), we obtain \(F_q\equiv 7\pmod{19}\), but Table 3 shows that \(7\) is not a Fibonacci residue modulo \(19\), contradiction.

If \(k=\!3\), then \(F_8=21\), hence \(n=23+F_q\). From \(n\equiv 0\pmod{23}\), we get \(F_q\equiv 0\pmod{23}\). Recalling from Table 1 that the Fibonacci sequence modulo \(23\) has period \(48\), and by an elementary computation we know that for \(0\le m\le 47\), \(23\mid F_m\) if and only if \(m=0,24\), we deduce \(24\mid q\), impossible for a prime.

If \(k\ge 4\), then \(16\mid 2^k\), so \(F_{2^k}\equiv 0\pmod 7\). Thus
\[
n\equiv 2+F_q\pmod 7.
\]
Since \(n\equiv 2\pmod 7\), we obtain \(F_q\equiv 0\pmod 7\). From Table 2, \(F_m\equiv 0\pmod 7\) if and only if \(8\mid m\). Hence \(8\mid q\), which is impossible for a prime \(q\).

All possibilities lead to contradictions. Therefore no \(n\in\mathcal{A}\) admits a representation of the form \(p+F_{2^k}+F_q\) with primes \(p,q\) and \(k\ge 0\). This proves the lemma.\qed
\end{proof}
\medskip
\noindent\textbf{Proof of Theorem~\ref{thm:pos}.}
\bpf
The assertion follows directly from Lemma~\ref{lem:arith}, which provides an
explicit infinite arithmetic progression of positive integers, none of which
can be expressed as \(p+F_{2^k}+F_q\) with primes \(p,q\) and \(k\ge0\).
The existence of such a progression immediately implies that the set of
integers not representable in this form has positive lower asymptotic density.\qed
\epf
\section{Proof of Theorem \ref{thm:pos}}
\begin{lem}\label{lem:main}
Let the Fibonacci sequence be defined by \(F_0=0,\;F_1=1\) and
\(F_{n+2}=F_{n+1}+F_n\) for \(n\ge 0\).
For non‑negative integers \(k_1,k_2\) and primes \(q_1,q_2\), the equation
\begin{equation}\label{eq:main-eq}
F_{2^{k_1}}+F_{q_1}=F_{2^{k_2}}+F_{q_2}
\end{equation}
has only the solutions
\[
(k_1,q_1)=(k_2,q_2),\qquad
(k_1,k_2,q_1,q_2)=(1,0,q,q),\qquad
(k_1,k_2,q_1,q_2)=(0,1,q,q),
\]
where \(q\) is an arbitrary prime.
\end{lem}

\begin{proof}
We may assume \(k_1\ge k_2\); the case \(k_2>k_1\) follows by symmetry.

If \(k_1=k_2\), then \eqref{eq:main-eq} gives \(F_{q_1}=F_{q_2}\), which forces
\(q_1=q_2\) since \(q_i\) are primes and \(F_n\) is strictly increasing for
\(n\ge2\). This yields the first family.

Suppose \(k_1>k_2\ge1\). Then \(2^{k_1}-2^{k_2}\) is even, and by
Lemma~\ref{lem:fib}(8) we obtain
\begin{equation}\label{eq:diff-k}
F_{2^{k_1}}-F_{2^{k_2}}
= F_{2^{k_1}-1}+F_{2^{k_1}-3}+\cdots+F_{2^{k_2}+1},
\end{equation}
which is a sum of non-consecutive odd-indexed Fibonacci numbers.
From \eqref{eq:main-eq} this equals \(F_{q_2}-F_{q_1}\).
Since the latter is positive, \(q_2>q_1\).

If \(q_1\) is an odd prime, then \(q_2\) is odd and \(q_2-q_1\) is even.
Applying Lemma~\ref{lem:fib}(8) again, we obtain
\[
F_{q_2}-F_{q_1}=F_{q_2-1}+F_{q_2-3}+\cdots+F_{q_1+1},
\]
which is a sum of non-consecutive even-indexed Fibonacci numbers.
This contradicts the uniqueness of the Zeckendorf representation
(Lemma~\ref{lem:fib}(5)).

If \(q_1=2\), then \(F_{q_1}=1=F_1\). Since \(q_2\) is an odd prime,
\(q_2-1\) is even; applying Lemma~\ref{lem:fib}(8) to \(F_{q_2}-F_1\),
we obtain
\[
F_{q_2}-F_{q_1}=F_{q_2}-F_1=F_{q_2-1}+F_{q_2-3}+\cdots+F_2,
\]
again a sum of non-consecutive even-indexed Fibonacci numbers,
contradicting \\Lemma~\ref{lem:fib}(5).
Thus no solutions occur when \(k_1>k_2\ge1\).

Now consider \(k_1>k_2=0\). Then \(F_{2^{k_2}}=F_1=1\), and \eqref{eq:main-eq} becomes
\begin{equation}\label{eq:main-3}
F_{2^{k_1}}-1=F_{q_2}-F_{q_1}.
\end{equation}
If \(k_1=1\), then \(F_{2^{k_1}}=1\), so \eqref{eq:main-3} gives
\(F_{q_1}=F_{q_2}\), hence \(q_1=q_2\), giving the second family.

Assume \(k_1\ge2\). Then \(2^{k_1}-2\) is even, so
\begin{equation}\label{eq:main-4}
F_{2^{k_1}}-1=F_{2^{k_1}}-F_2
= F_{2^{k_1}-1}+F_{2^{k_1}-3}+\cdots+F_3,
\end{equation}
which is a sum of non-consecutive odd-indexed Fibonacci numbers
(all indices at least~3).
If \(q_1=2\), then \eqref{eq:main-3} gives \(F_{2^{k_1}}=F_{q_2}\), which is
impossible because \(2^{k_1}\) is composite for \(k_1\ge2\) while \(q_2\) is prime.
If \(q_1\) is odd, then \(q_2>q_1\) and both are odd, so
\[
F_{q_2}-F_{q_1}=F_{q_2-1}+F_{q_2-3}+\cdots+F_{q_1+1},
\]
a sum of even-indexed numbers, again contradicting Lemma~\ref{lem:fib}(5).
Hence no further solutions arise when \(k_2=0\) and \(k_1\ge2\).

Finally, if \(k_2>k_1\), interchanging the indices yields the symmetric
counterpart of the second family, namely \(k_1=0,\;k_2=1,\;q_1=q_2\),
which is the third family. This completes the proof.\qed
\end{proof}

\begin{lem}\label{lem:Rx}
Let \(r(n)=\#\{(p,k,q)\in\mathcal{P}\times\mathbb{N}_0\times\mathcal{P}:
n=p+F_{2^k}+F_q\}\) and
\(R(x)=\sum\limits_{n\le x}r(n)\).
Then \(R(x)\sim x\) as \(x\to\infty\).
\end{lem}

\begin{proof}
Denote by \(\pi(x)\) the number of primes \(p\le x\).
The Prime Number Theorem gives \(\pi(x)\sim x/\log x\).
From Lemma~\ref{lem:fib}(4) we have \(L(x)\sim\log x\) and \(W(x)\sim\log x\)
for large \(x\).
Notice \(R(x)\) counts triples \((p,k,q)\) with \(p,q\) prime, \(k\ge0\),
and \(p+F_{2^k}+F_q\le x\).

On the one hand, if \(p+F_{2^k}+F_q\le x\), then \(F_{2^k},F_q,p\le x\).
Combining this with Lemma~\ref{lem:fib}(2) and (4), we obtain
\[
R(x)\le\bigl(\log_2 L(x)+1\bigr)\cdot\pi(L(x))\cdot\pi(x)
\ll(\log\log x)\cdot\frac{\log x}{\log\log x}\cdot\frac{x}{\log x}=x.
\]

On the other hand, the conditions
\(p\le x/3\), \(F_{2^k}\le x/3\), and \(F_q\le x/3\)
imply \(p+F_{2^k}+F_q\le x\).
By Lemma~\ref{lem:fib}(3), \(F_m\le x/3\) holds whenever \(m\le W(x/3)\).
Therefore,
\[
\begin{aligned}
\#\{k:F_{2^k}\le x/3\}
&\ge\#\{k:2^k\le W(x/3)\}\gg\log W(x/3)\sim\log\log x,\\
\#\{q\text{ prime}:F_q\le x/3\}
&\ge\pi\bigl(W(x/3)\bigr)\sim\frac{\log x}{\log\log x}.
\end{aligned}
\]
Together with \(\pi(x/3)\sim x/\log x\) we obtain
\[
R(x)\gg(\log\log x)\cdot\frac{\log x}{\log\log x}\cdot\frac{x}{\log x}=x.
\]
Hence \(R(x)\sim x\).\qed
\end{proof}
\begin{lem}\label{lem:second-moment}
$\displaystyle\sum_{n\le x}r(n)^2\ll x$.  
\end{lem}

\bpf
It is obvious that
\begin{equation}\label{eq:r2-count}
\begin{aligned}[b]
\sum_{1 \leq n \leq x} r(n)^2
&= \#\bigl\{(p_1, k_1, q_1, p_2, k_2, q_2) \mid p_1+F_{2^{ k_1}}+F_{q_1}=p_2+F_{2^{k_2}}+F_{q_2} \leq x\bigr\}.
\end{aligned}
\end{equation}
Let \(
h := F_{2^{k_1}} + F_{q_1} - F_{2^{k_2}} - F_{q_2},
\) then \(p_2 = p_1 + h\).
We now classify the solutions \((p_1, k_1, q_1, p_2, k_2, q_2)\) according to the value of
\(h = F_{2^{k_1}} + F_{q_1} - F_{2^{k_2}} - F_{q_2}\) into the following four cases:
\[
\text{Case I: } h = 0;\quad
\text{Case II: } h < 0;\quad
\text{Case III: } h > 0,\text{ odd};\quad
\text{Case IV: } h > 0,\text{ even}.\]
\textbf{Case I:} If $h = 0$, then $F_{2 ^{k_1}}+F_{q_1}=F_{2^{k_2}}+F_{q_2}$.
By Lemma~\ref{lem:main}, this yields $(k_1, q_1) = (k_2, q_2)$; or $k_1 = 0$, $k_2 = 1$, $q_1 = q_2$; or $k_1 = 1$, $k_2 = 0$, $q_1 = q_2$.
In this case, for fixed $p_1, k_1, q_1$, the condition $h = 0$ determines $k_2$ and $q_2$, and then the equality
$p_2 = p_1 - F_{2^k} - F_{q_1} - F_{2^{k 2}} - F_{q_2}$
determines $p_2$.  Therefore, in this situation, the number of solutions
$(p_1, k_1, q_1, p_2, k_2, q_2)$ satisfying
$F_{2^{k_1}} + F_{q_1} - F_{2^{k_2}} - F_{q_2} =0$ and $p_1 + F_{2^{ k_1}} + F_{q_1} = p_2 + F_{2^{k_2}} + F_{q_2} \le x$
is given by
\[
\sum_{1 \leqslant n \leqslant x} r(n) \ll x .
\]\\ \textbf{Case II:} If \(h<0\), then \(-h = F_{2^{k_2}} + F_{q_2} - F_{2^{k_1}} - F_{q_1}\) and \(-h = p_1 - p_2\); by symmetry, this has the same form as the case \(h>0\) after swapping the indices. Hence, without loss of generality, we assume \(h > 0\).

Since \(3\nmid 2^k\) and the Fibonacci sequence modulo \(2\) has period \(3\), it follows that \(F_{2^k}\) is odd for every \(k\ge 0\) and $F_{q_1},F_{q_2}$ is odd for every $3\nmid q_1, 3\nmid q_2.$ \\
\textbf{Case III:} If $h>0$ is odd, then from $h = p_2 - p_1$ we deduce that $p_1 = 2$ and ( $q_1=3$ or $q_2=3$).
Therefore in this case, from
\[
2 + F_{2^{k_{1}}} + F_{q_1} = p_2 + F_{2^{k_2}} + F_{q_2}
\]
we see that once $k_1, k_2,  p_2$ are determined, $q_1,q_2$ is determined as well. 
Therefore, the number of solutions \((p_1, k_1, q_1, p_2, k_2, q_2)\) such that
\(F_{2^{k_1}} + F_{q_1} - F_{2^{k_2}} - F_{q_2}\) is odd and
\(p_1 + F_{2^{k_1}} + F_{q_1} = p_2 + F_{2^{k_2}} + F_{q_2} \le x\)
is at most
\[
\frac{\log L(x)}{\log 2} \cdot \frac{\log L(x)}{\log 2} \cdot \frac{L(x)}{\log L(x)} = o(x).
\]\\ \textbf{Case IV:} Therefore, it now suffices to estimate the contribution from the case when \(h\) is positive and even.
Since \(h = p_2 - p_1\), we apply Lemma~\ref{lem:BD} to bound the number of solutions
\((p_1, k_1, q_1, p_2, k_2, q_2)\) to
\[
p_1 + F_{2^{k_1}} + F_{q_1} = p_2 + F_{2^{k_2}} + F_{q_2} \le x,
\]
where \(h = F_{2^{k_1}} + F_{q_1} - F_{2^{k_2}} - F_{q_2}\) is a positive even integer.
Thus
\begin{equation}\label{eq:bd-applied}
\begin{aligned}[b]
&\#\{(p_1,k_1,q_1,p_2,k_2,q_2): p_1+F_{2^{k_1}}+F_{q_1}=p_2+F_{2^{k_2}}+F_{q_2}\le x,h>0 \;\text{even};\}\\
&\ll \sum_{
\substack{
(k_1, q_1, k_2, q_2) \\[2pt]
F_{2^{k_1}},\, F_{2^{k_2}},\, F_{q_1},\, F_{q_2} \le x \\[2pt]
F_{2^{k_1}} + F_{q_1} - F_{2^{k_2}} - F_{q_2} > 0,\ \text{even}
}}
\frac{x}{\log^2 x}
\prod_{\substack{p \text{ prime} \\ p \mid F_{2^{k_1}} + F_{q_1} - F_{2^{k_2}} - F_{q_2}}}
\!\left(1 + \frac{1}{p}\right) \\
&\;\text{For brevity, denote the above} \;= \frac{x}{\log^2 x} \sideset{}{'}\sum_{(k_1, q_1, k_2, q_2)} \prod_{p \mid h}\!\left(1 + \frac{1}{p}\right).
\end{aligned}
\end{equation}

Since $h = F_{2^{k_1}} + F_{q_1} - F_{2^{k_2}} - F_{q_2} \le x$, let $s$ denote the number of prime factors of $h$ (counted without multiplicity) that are at least $\log x$.
Because each of these prime factors is at least $\log x$, we have $h \ge (\log x)^s$, and thus $x \ge (\log x)^s$.
It follows that $s \le \frac{\log x}{\log\log x}$.
Hence
\[
\prod_{\substack{p \mid h \\ p \ge \log x}} \left(1+\frac{1}{p}\right) 
\le \left(1+\frac{1}{\log x}\right)^{\frac{\log x}{\log\log x}} 
\to 1 \quad (x \to +\infty).
\]
Thus it suffices to estimate
\[
\sideset{}{'}\sum_{(k_1, q_1, k_2, q_2)}
\prod_{\substack{p \mid h \\ p \le \log x}} \Bigl(1 + \frac{1}{p}\Bigr) \ll \log^2 x.
\]
We will use $P(d)$ to denote the largest prime factor of $d$.\\
Consequently,
\begin{equation}\label{eq:swap-sum}
\begin{aligned}[b]
\sideset{}{'}\sum_{(k_1, q_1, k_2, q_2)}
\prod_{\substack{p \mid h \\ p \le \log x}} \Bigl(1 + \frac{1}{p}\Bigr)
= \sideset{}{'}\sum_{(k_1, q_1, k_2, q_2)}
\sum_{\substack{d \mid h \\ P(d) \le \log x}} \frac{\mu^2(d)}{d}.
\end{aligned}
\end{equation}
Interchanging the order of summation, we obtain
\[
\sum_{\substack{d \geq 1 \\ P(d) \le \log x}} \frac{\mu^2(d)}{d}  \; \# S_{d,x},
\]
where
\[
S_{d,x} = \left\{ (k_1, q_1, k_2, q_2) : 
\begin{array}{l}
F_{2^{k_1}}, F_{2^{k_2}}, F_{q_1}, F_{q_2} \le x,\\[6pt]
F_{2^{k_1}} + F_{q_1} - F_{2^{k_2}} - F_{q_2} > 0\\[6pt]
\text{even, and divisible by } d
\end{array}
\right\}.
\]
We define
\[
H_{d,x} = \left\{
(k_1, q_1, k_2, q_2) : 
\begin{aligned}
&F_{2^{k_1}}, F_{2^{k_2}}, F_{q_1}, F_{q_2} \le x,\\
&F_{2^{k_1}}+F_{q_1}-F_{2^{k_2}}-F_{q_2} \equiv 0 \pmod d
\end{aligned}
\right\}.
\]
Clearly $\# S_{d,x} \le \# H_{d,x}$.

By Lemma~\ref{lem:fib}(6), the Fibonacci sequence is purely periodic modulo any integer.
We introduce the following notation.
Let \(k(d)\) denote the minimal positive period of the Fibonacci sequence modulo \(d\).
For a prime \(p\), we 
define
\[
v(d) = \max\{\,k(p) : p\mid d,\ p\ \text{prime}\,\}
\]
and
\[
LP(d) = \{\,p\mid d : p\ \text{prime and}\ k(p) \ge k(p')\ \text{for all primes}\ p'\mid d\,\}.
\]
For each positive integer $d$, we fix a prime $p_d \in LP(d)$.
Now fix \(k_1, k_2, q_2\) with \(F_{2^{k_1}}, F_{2^{k_2}}, F_{q_2} \le x\), and consider the congruences
\begin{align}
F_{q_1} &\equiv F_{2^{k_2}} + F_{q_2} - F_{2^{k_1}} \pmod d,
\quad F_{q_1} \le x,
\label{eq:mod-d}
\\
F_{q_1} &\equiv F_{2^{k_2}} + F_{q_2} - F_{2^{k_1}} \pmod{p_d},
\quad F_{q_1} \le x.
\label{eq:mod-pd}
\end{align}
The number of $q_1$ satisfying \eqref{eq:mod-d} is at most the number of solutions of \eqref{eq:mod-pd}.

By Lemma~\ref{lem:fib} we have $2\log x < L(x) < 4\log x$ for all $x\ge 3$,
and $F_n\le x$ implies $n\le L(x)$.
Set $l = F_{2^{k_2}} + F_{q_2} - F_{2^{k_1}}$; in particular, $l$ is then fixed as well.
In this case, consider the congruence
\[
F_{q_1} \equiv l \pmod{p_d},\quad
q_1 \le L(x).
\]
Since $p_d\in LP(d)$, the minimal positive period of $F_n$ modulo the prime $p_d$ is $v(d)$.
Moreover, by Lemma~\ref{lem:fib}(7), every residue occurs at most four times in one full period of the Fibonacci numbers modulo any prime.
Consequently, the indices $q_1$ satisfying the congruence $F_{q_1} \equiv l \pmod{p_d}$ belong to the union of at most four residue classes modulo $v(d)$. 
\begin{enumerate}[label=(\roman*)]
\item Suppose $v(d) < \log x$. Then, because $L(x) > 2\log x$, we clearly have $v(d) < L(x)$.  
Consider an arbitrary residue class modulo $v(d)$.
\begin{itemize}
    \item[(a)] If the representative is not coprime to $v(d)$, the class contains at most one prime.
    \item[(b)] If the representative is coprime to $v(d)$, then by Lemma~\ref{lem:BT} the number of primes $\le L(x)$ (note that $L(x)>v(d)$) in this class is bounded by
    \[
    \frac{2 L(x)}{\varphi(v(d)) \log\bigl( L(x) / v(d) \bigr)} .
    \]
\end{itemize}
From the earlier discussion, the indices $q_1$ with $F_{q_1} \equiv l \pmod{p_d}$ (for some fixed $p_d\in LP(d)$) lie in the union of at most four residue classes modulo $v(d)$.  
Therefore, when $v(d) < \log x$, we obtain
\begin{equation}\label{eq:q1-bound-small}
\small
\#\{\, q_1 : F_{q_1} \equiv l \pmod{p_d},\ q_1\le L(x) \,\}
\le 4 \max\{ 1,\ \frac{2 L(x)}{\varphi(v(d)) \log(L(x)/v(d))} \}.
\end{equation}
\item Suppose $v(d) \ge \log x$.
Since $F_{q_1} \le x$, we have $q_1 \le L(x) < 4\log x \le 4v(d)$.
As $q_1$ is a prime, $q_1 \ge 2$, so $2 \le q_1 \le 4v(d)$.
By Lemma~\ref{lem:fib}(7) (each residue occurs at most four times in one period),
the congruence $F_{q_1} \equiv l \pmod{p_d}$ has at most $4 \times 4 = 16$
solutions in $q_1$ within this range.
\end{enumerate}
Observe that the total number of possible triples $(k_1, k_2, q_2)$ is at most
\[
\frac{\log L(x)}{\log 2} \cdot \frac{\log L(x)}{\log 2} \cdot \frac{L(x)}{\log L(x)} .
\]
Therefore,
\begin{equation}\label{eq:Hd-bound}
\small
\# H_{d,x}
\le \frac{L(x) \log L(x)}{(\log 2)^2}
\times
\begin{cases}
16, & v(d) \ge \log x,\\[4pt]
4 \max\{ 1,\ \frac{2 L(x)}{\varphi(v(d)) \log(L(x)/v(d))} \}, & v(d) < \log x.
\end{cases}
\end{equation}

Inserting the bound from \eqref{eq:Hd-bound} into the sum, we obtain 
\begin{equation}\label{eq:main-sum}
\small
\sum_{\mathclap{\substack{d \ge 1 \\ P(d) \le \log x}}}
\frac{\mu^2(d)}{d} \, \# H_{d,x}
\le
\frac{L(x) \log L(x)}{(\log 2)^2}
\Biggl(
16 \sum_{\mathclap{\substack{1 \le d \le x \\ P(d) \le \log x \\ v(d) \ge \log x}}}
\frac{\mu^2(d)}{d}
+ 4 \sum_{\mathclap{\substack{1 \le d \le x \\ P(d) \le \log x \\ v(d) < \log x}}}
\frac{\mu^2(d)}{d}
\max\!\bigl\{1,\; \frac{2 L(x)}{\varphi(v(d)) \log(L(x)/v(d))}\bigr\}
\Biggr).
\end{equation}Applying the elementary bound
\(\max\{1,\, A\} \le 1 + A\)
with \(A = \frac{2 L(x)}{\varphi(v(d)) \log(L(x)/v(d))}\)
and using the estimate for \(\varphi (v(d))\) from Lemma~\ref{lem:Hardy},
we obtain
\[
\small
\sum_{\mathclap{\substack{d \ge 1 \\ P(d) \le \log x}}} \frac{\mu^2(d)}{d} \, \# H_{d,x}
\le \frac{L(x) \log L(x)}{(\log 2)^2}
\Biggl(
20 \sum_{\mathclap{\substack{1 \le d \le x \\ P(d) \le \log x}}} \frac{\mu^2(d)}{d}
+ 8 \sum_{\mathclap{\substack{1 \le d \le x \\ P(d) \le \log x \\ v(d) < \log x}}} 
\frac{\mu^2(d)L(x) \log\log v(d)}{d \, v(d)\log(L(x)/v(d))}
\Biggr).
\]

By Mertens' theorem,
\[
\sum_{\substack{1 \le d \le x \\ P(d) \le \log x}} \frac{\mu^2(d)}{d}
= \prod_{p \le \log x} \Bigl(1 + \frac{1}{p}\Bigr) \le \log\log x .
\]

It remains to handle the sum
\[
S := \sum_{\substack{1 \le d \le x \\ P(d) \le \log x \\ v(d) < \log x}} \frac{\mu^2(d) \cdot \log\log v(d)}{d} \cdot \frac{L(x)/v(d) }{\log(L(x)/v(d))}.
\]
We split the range of $v(d)$ with two fixed constants $0<\varepsilon<1$ and $0<\eta<\frac34$.

\noindent\textbf{Case 1:} $(\log x)^\varepsilon < v(d) < \log x$.  
Then
\[
\frac{L(x)}{v(d)} = \frac{L(x)}{v(d)^\eta \, v(d)^{1-\eta}}
\le \frac{4 \log x}{(\log x)^{\varepsilon\eta} \, v(d)^{1-\eta}}
= \frac{4(\log x)^{1-\varepsilon\eta}} {v(d)^{1-\eta}} .
\]
Note that the sum \(S\) we are bounding is taken only over those \(d\) with
\(v(d) < \log x\) and \(P(d) \le \log x\).
Since \(L(x) > 2\log x\) for all \(x \ge 3\), it follows that
\[
\frac{L(x)}{v(d)} > \frac{2\log x}{v(d)} > \frac{2\log x}{\log x} = 2.
\]
In this case, we have
\[
\frac{L(x)}{v(d)\log(L(x)/v(d))}
< \frac{4(\log x)^{1-\varepsilon\eta}}{v(d)^{1-\eta}\log 2}.
\]
\noindent\textbf{Case 2:} $v(d) \le (\log x)^\varepsilon$.  
Here
\[
\frac{L(x)}{v(d)} > \frac{2 \log x}{(\log x)^\varepsilon} = 2(\log x)^{1-\varepsilon},
\]
so $\log\frac{L(x)}{v(d)} > \log 2 + (1-\varepsilon)\log\log x$, and consequently
\[
\frac{1}{\log\frac{L(x)}{v(d)}} \ll \frac{1}{(1-\varepsilon)\log\log x}.
\]
Thus
\[
\frac{L(x)}{v(d) \log(L(x)/v(d))} \ll \frac{4 \log x}{(1-\varepsilon) v(d) \log\log x}.
\]

Combining the estimates obtained in Case~1 and Case~2 above, we obtain
\[
S \ll 
\frac{4(\log x)^{1-\varepsilon\eta}}{\log 2}
\sum_{\mathclap{\substack{1 \le d \le x \\ P(d) \le \log x \\ (\log x)^\varepsilon < v(d) < \log x}}}
\frac{\mu^2(d) \log\log v(d)}{d \, v(d)^{1-\eta}}
+ \frac{4 \log x}{(1-\varepsilon) \log\log x}
\sum_{\mathclap{\substack{1 \le d \le x \\ P(d) \le \log x \\ v(d) \le (\log x)^\varepsilon}}}
\frac{\mu^2(d) \log\log v(d)}{d \, v(d)} .
\]

Since $v(d) \ge 3$ (a Fibonacci period modulo a prime cannot be $1$ or $2$) and $\log\log x \le x^{1/4}$ for all $x \ge 3$, we have
\[
\frac{\mu^2(d) \log\log v(d)}{d \, v(d)} \le \frac{\mu^2(d)}{d \, v(d)^{\frac{3}{4}}}, \qquad
\frac{\mu^2(d) \log\log v(d)}{d \, v(d)^{1-\eta}} \le \frac{\mu^2(d)}{d \, v(d)^{\frac34 - \eta}} .
\]

By Lemma~\ref{lem:chen}, the series $\sum\limits_{d=1}^\infty \frac{\mu^2(d)}{v(d)^\alpha}$ converges for every $\alpha > 0$.  
Choosing $\alpha_1 = \frac34 - \eta > 0$ and $\alpha_2 = \frac34 > 0$, we obtain absolute constants $M_1, M_2 > 0$ such that
\[
\sum_{\substack{1 \le d \le x \\ P(d) \le \log x \\ (\log x)^\varepsilon < v(d) < \log x}} \frac{\mu^2(d) \log\log v(d)}{d \, v(d)^{1-\eta}} \le M_1,
\qquad
\sum_{\substack{1 \le d \le x \\ P(d) \le \log x \\ v(d) \le (\log x)^\varepsilon}} \frac{\mu^2(d) \log\log v(d)}{d \, v(d)} \le M_2 .
\]

Therefore,
\[
S \ll M_{1}(\log x)^{1-\varepsilon\eta} + M_{2}\frac{\log x}{\log\log x}.
\]

Collecting all the estimates and recalling $L(x) < 4\log x$, we conclude
\[
\sum_{\substack{1 \le d \le x \\ P(d) \le \log x}} \frac{\mu^2(d)}{d} \, \# H_{d,x}
\le \frac{L(x) \log L(x)}{(\log 2)^2}
\Bigl( 20 \log\log x + 8 S \Bigr)
\ll \log^2 x .
\]

Combining all four cases, we obtain
\[
\sum_{1 \le n \le x} r(n)^2 \ll x .
\]\qed
\epf

\medskip
\noindent\textbf{Proof of Theorem~\ref{thm:pos}.}
\begin{proof}
By Lemma~\ref{lem:Rx} we have $\sum\limits_{n\le x} r(n) \gg x$, and by
Lemma~\ref{lem:second-moment} we have $\sum\limits_{n\le x} r(n)^2 \ll x$.
Applying the Cauchy–Schwarz inequality,
\[
\sum_{n\le x} r(n) \le \Bigl(\sum_{\substack{n\le x \\ r(n)\ge 1}} 1\Bigr)^{\!1/2}
\Bigl(\sum_{n\le x} r(n)^2\Bigr)^{\!1/2},
\]
the two estimates above force
\[
\#\{n\le x : r(n)\ge 1\} \gg x,
\]
i.e., the set of integers representable as $p+F_{2^k}+F_q$ has positive lower
asymptotic density.\qed
\end{proof}

\end{document}